\documentclass[a4paper, 12pt]{article}
\usepackage{enumerate, theorem}
\usepackage{amsmath, amsfonts, amssymb}
\usepackage[height=22.5cm, width=15cm]{geometry}
\usepackage[all]{xy}
\usepackage{mathrsfs, yfonts}
\usepackage{graphicx}

 \def\C{{\mathbb C}}

\def\bP{{\mathbb P}}

\newcommand{\parag}[1]{\paragraph{\sc{#1.}} }

\DeclareMathOperator{\I}{\mathcal{I}}
\DeclareMathOperator{\J}{\mathcal{J}}

\newtheorem{thm}{Theorem}[subsection]
\newtheorem{defn}[thm]{Definition}
\newtheorem{cor}[thm]{Corollaire}
\newtheorem{prop}[thm]{Proposition}
\newtheorem{lemma}[thm]{Lemma}

\begin{document}

\date{26/03/15}

\author{Daniel Barlet\footnote{Invited professor to the KIAS (Seoul), Institut Elie Cartan : Alg\`{e}bre et G\'eom\`{e}trie,  \newline
Universit\'e de Lorraine, CNRS UMR 7502  and  Institut Universitaire de France.}.}

\title{Some examples due to H. Hironaka}

\maketitle

\section*{Abstract}

The aim of this paper is to explain the construction by H. Hironaka [H.61] of a holomorphic (in fact ''algebraic'') family of compact complex manifolds  parametrized by $\C$ such for all $s \in \C\setminus \{0\}$ the fiber is projective, but such that the fiber at the origin in non k{\"a}hlerian, to mathematicians which are not algebraic geometers. We also explain why it is not possible to make in the same way such  an example with fiber at $0$ a simpler example of non k{\"a}hlerian Moishezon manifold which is also due to H. Hironaka (see section 5).

\parag{AMS CLASSIFICATION 2010} 32-G-10, 32-J-27, 14-C-25, 14-D-06.

\parag{Key words} Family of compact complex manifolds - Non k\"ahlerian Moishezon manifold.

\tableofcontents

\section*{Introduction}

The main goal of this work is to provide an introduction to the article of H. Hironaka [H.61]. 
It concerns the construction of a holomorphic (in fact ``algebraic'') family of compact complex manifolds  parametrized by $\C$ with the following property. For all $u \in \C\setminus \{0\}$ the fiber is projective, but the fiber at the origin in not K{\"a}hler (it is however a Moishezon manifold). This famous example serves to show that large deformations of projective manifolds may not be K\"ahler.
\\

The central fiber of the said family is quite complicated and one may wonder if it is not possible to 
simplify the construction, so as to understand better the reason why the K\"ahler property may disappear in large deformations. A possible candidate for the central fiber 
is the simpler example of non k{\"a}hlerian Moishezon manifold, 
also due to H. Hironaka, see section 5. 
We will explain here the reason why it is not possible to obtain, with an analoguous construction, a family whose fiber at $0$ is this precise 
simpler example. 
\\

It was not our intention to present here a complete proof of Hironaka's construction. Instead, we have rather 
detailed some particular points of the original article which --in our opinion-- are not ``explicit'' enough for a ``generic'' reader, although they are well known to specialists. In the same spirit, of helping the reader to gain some intuition about the geometry of the family 
obtained, we have tried to draw a few two-dimensional pictures in the appendix (revealing many-dimensional phenomenons...).
\\

This article is organized as follows: after giving a rough idea of Hironaka's construction in section 1 we recall some basic properties of the blow-up of a coherent ideal which will be useful later on. In section 3 we will prove a technical result, which plays a key role in a glueing procedure used in the construction. The properties of the central fiber of the family are further analyzed in section 4: we show e.g. that this variety is smooth -- it is one of the delicate points of [H.61], and that it is not K\"ahler.
We remark that there are simpler examples of non-K\"ahler manifolds than the central fiber mentioned above, and thus it may be tempting to try use them in order to obtain non-projective fibers as large deformations of 
projective manifolds.
In section 5 we explain why it is not possible to do so, at least for the class of non-K\"ahler 3-folds
we consider here. \\
We conclude by exhibiting in section 6 a non-K\"ahler Moishezon manifold containing two disjoint smooth curves such that the blow-up of any of them is projective. This example shows that the fact that this manifold is not K\"ahler cannot be ``localized''.\\

I want to thank Mihai Paun for his invitation to the KIAS and also for helping me to improve this text.

\section{The construction: a quick overview}

The initial point of H. Hironaka's article is to consider in the $3-$dimensional complex projective space $\mathbb{P}_{3}$ the three curves given in the homogeneous coordinates $Y_{0}, Y_{1}, Y_{2},Y_{3}$ by the following equations (see picture 4)
\begin{align*}
& F_{1} := \{ Y_{2} = Y_{3} = 0 \} \quad ( = C_{1} \quad {\rm  on  \  pictures}) \\
& F_{2 } := \{ Y_{0}.Y_{1} + Y_{1}.Y_{2} + Y_{2}.Y_{0} = Y_{3} = 0 \} \quad   ( = C_{2} )\\
&  F_{3}(u)  := \{ (Y_{1} + Y_{3}).(Y_{0} + u.Y_{1}) + Y_{1}.Y_{3} = Y_{2} = 0 \} \quad  (  F_{3}(0) = C_{3} )
\end{align*}
where $u \in \C$ is a parameter.\\
 Note that for $u = 0$  these three smooth  curves (a line and two conics) has two common points $P := (1,0,0,0)$ and $Q := (0,1,0,0)$. For $u \not= 0$ the conic $F_{3}(u)$ still contains $P$ but not $Q$. In fact $P$ is the only common point between $F_{3}(u), F_{1}$ and  $F_{2}$ for $u \not= 0$, as $F_{3}(u) \cap F_{1} = \{P, Q(u)\} $ where $Q(u) := (-u, 1, 0, 0) \not\in F_{2}$.\\

 We first blow-up the ideal corresponding to the union $C:= C_1\cup C_2\cup C_3$ in $\bP_3$. The central fiber of the family we will construct is obtained by a glueing procedure starting from the resulting (projective) variety, as follows.
 
 Let
 $$U_{0} := \{ Y_{0}\not= 0\}\subset \bP_3$$
 be an affine open set, 
 and let $(x_1,\dots x_3)$ be the corresponding inhomogeneous coordinates. The next step is to blow-up the
 ideal 
 $$\mathcal{H} := (x_{1}.x_{2}, x_{2}.x_{3}, x_{1}.x_{3}).(x_{1}, x_{2}.x_{3}).(x_{2}, x_{1}.x_{3})^{2}.(x_{3}, x_{1}.x_{2})^{2}$$
 in $U_0$. Let $U_1$ be the affine set 
$\{ Y_{1} \not= 0\}$; it turns out that on the intersection $U_0\cap U_1$, 
these two operations agree. As a consequence, we can modify the first (globally defined) blow-up 
of $\bP_3$ along $C$ in such a 
way that the resulting manifold denoted by $V_0$ will coincide with 
the second blow-up in $U_0$.
We will show that $V_{0}$ is not K{\"a}hler, so in particular it is not projective. This is proved by showing that there exists  an effective curve in $V_{0}$ which is homologous to zero (cf. the end of section 3). We remark, however, that $V_0$ is a modification of the projective space 
$\mathbb{P}_{3}$.
\smallskip

Next we will perform a similar construction with the parameter $u \in \C$ in order to produce an holomorphic family of compact complex manifolds over $\C$ such that the fibers over $u \not= 0$ are projective but the fiber over $0$ (which is the previously constructed $V_{0}$)  is not K{\"a}hler.
 This second step is performed by using basically the same procedure, simply by adding the parameter $u \in \C$. On $\C\times U_{0}$ we will replace the inhomogeneous coordinates $x_{1}, x_{2}, x_{3}$ on $U_{0}$ by
  $$\tilde{x}_{1} := (1 + u.x_{1})(x_{1} + x_{2}+ x_{3} + x_{1}.x_{2}) + x_{1}.x_{3}, \quad  \tilde{x}_{2} := x_{2}, \quad   \tilde{x}_{3} := x_{3}$$
and then follow the same line of arguments as in the absolute case: \\
 the map $(x_{1}, x_{2}, x_{3}, u) \mapsto (\tilde{x}_{1},  \tilde{x}_{2},  \tilde{x}_{3}, u) $ is a local holomorphic change of variables near the origin in $\C^{3}\times \C = U_{0}\times \C$, so the blow-up of the ideal $\tilde{\mathcal{H}}$ obtained from $\mathcal{H}$ by changing $x_{1}, x_{2}, x_{3}$ in $\tilde{x}_{1},  \tilde{x}_{2},  \tilde{x}_{3}$ gives a family parametrized by $u \in \C$ of smooth complex manifolds with fiber $U_{0}\cap V_{0}$ for $u = 0$.\\
We check in details the smoothness of $U_{0}\cap V_{0}$ in the section 4.\\
The projectivity of the fibers over $u \not= 0$ is  deduced from the fact that the corresponding  manifold (built after glueing) can be obtained directly by blowing up a coherent ideal in $ \mathbb{P}_{3} \times \{u\}$.\\
 We give an example in section 5 where an analogous construction ``in family'' do not give the expected fiber at $u = 0$ and this explains why three curves are needed through the point $P$ in Hironaka's example.\\

\section{Blowing up coherent ideals}

In the present text, an ``ambiant '' complex space is always assumed to be  reduced. But it will be important to consider the ``natural'' complex structure on a fiber of a morphism (between reduced complex spaces) which is not always reduced. Recall that if $f : Y \to X$ is a holomorphic map between (reduced) complex spaces, the fiber of $f$ at $x \in X$ is the complex subspace of $Y$ associated to the coherent ideal $f_{str}^{*}(\frak{M}_{x})$\footnote{See the definition below.} of $\mathcal{O}_{Y}$ where $\frak{M}_{x}$ is the maximal ideal of $\mathcal{O}_{X}$  at $x$.\\

Also a modification of a complex space will be always assumed to be proper.\\

When $f  : Y \to X$ is a morphism of complex spaces and $\mathcal{I}$ a sheaf of ideals in $\mathcal{O}_{X}$; the pull-back by $f$ of the  ideal sheaf $\mathcal{I}$, denoted by $f_{st}^{*}(\mathcal{I})$, is the image in $\mathcal{O}_{Y}$ of the ``usual'' pull-back sheaf $f^{*}(\mathcal{I}) := f^{-1}(\mathcal{I}) \otimes_{f^{-1}(\mathcal{O}_{X})} \mathcal{O}_{Y}$. If $\mathcal{I}$ is locally generated on an open set $U$ by holomorphic functions $g_{1}, \dots, g_{k}$ then the ideal $f_{st}^{*}(\mathcal{I})$ is generated on $f^{-1}(U)$ by the holomorphic functions $g_{1}\circ f, \dots, g_{k}\circ f$.\\

To begin, recall that if $X$ is a complex space and $\mathcal{I}$ a coherent ideal sheaf in $\mathcal{O}_{X}$ such that $V(\mathcal{I})$, the common zero set of sections of $\mathcal{I}$\footnote{Or the support of $\mathcal{O}_{X} /\mathcal{I}$.},  has no interior point in $X$, the blow-up  of $\mathcal{I}$ is, by definition, the (proper) modification $\tau : \tilde{X} \to X$ such that the ideal $\tau_{st}^{*}(\mathcal{I})$ is locally principal and  which satisfies the following universal property:
\begin{itemize}
\item For any (proper) modification $\theta : Y \to X$ such that the ideal $\theta_{st}^{*}(\mathcal{I})$ is a locally principal ideal sheaf, there exists an unique holomorphic map $\sigma : Y \to \tilde{X}$ such that $ \tau\circ \sigma = \theta$.
\end{itemize}

For any coherent ideal $\mathcal{I}$ in $\mathcal{O}_{X}$ such that $V(\mathcal{I})$ has no interior point in $X$,  the blow-up exists and is unique. It can be constructed as follows :\\

Assume that $\mathcal{I}$ is locally generated on an open set $U \subset X$ by holomorphic functions $g_{1}, \dots, g_{k}$ ; then, if $Z$ is the analytic subset of $U$ defined by $ \{g_{1}= \cdots = g_{k} = 0 \}$, we have a holomorphic map $ G : U \setminus Z \to \mathbb{P}_{k-1}$ defined in homogeneous coordinates by $x \mapsto (g_{1}(x), \dots, g_{k}(x))$. Then the closure of the graph $ \Gamma_{0}$  of $G$ in $U \times \mathbb{P}_{k-1}$ is an analytic subset and its projection on $U$ is a proper modification of $U$. The easiest way to prove these assertions is, assuming $U$ irreducible for simplicity, to consider the irreducible component $\Gamma$  of the closed analytic subset
$$ Y := \{(x, z) \in U \times \mathbb{P}_{k-1} \ / \  rk(G(x), z) = 1 \}$$
which contains the graph of $G$. Then it is easy to see that it coincides with the closure of $\Gamma_{0}$. It is clearly proper over $U$ and the pull-back of $\mathcal{I} := (g_{1}, \dots, g_{k})$ on $\Gamma$ is locally principal : if $z_{1}, \dots, z_{k}$ are the homogeneous coordinates in $\mathbb{P}_{k-1}$, let $\Omega_{i} := \{z_{i}\not= 0\} \subset \mathbb{P}_{k-1}$. Then on $(U \times \Omega_{i})\cap \Gamma$ we have $g_{j} = \frac{z_{j}}{z_{i}}.g_{i}$ for each $j \in [1,k]$. So $\I$ is generated by $g_{i}$ on $(U \times \Omega_{i})\cap \Gamma$.

\parag{Remarks}\begin{enumerate}
\item For $\J$ a coherent ideal sheaf  consider $ \mathcal{I} = \varphi.\mathcal{J}$ where $\varphi$ is a holomorphic function on $X$, which is not a zero divisor; then there is a canonical isomorphism of modifications of $X$ between the blow-up of $\mathcal{I}$ and the blow-up of $\mathcal{J}$.
\item If $X$ is a complex manifold and $\mathcal{I}$ is the (reduced) ideal sheaf of a closed complex submanifold $V$ with no interior point in $X$, then the blow-up $\tilde{X}$  of $\mathcal{I}$ is a complex manifold and the map $\tau : \tau^{-1}(V) \to V$ is the projection of the projectivized normal bundle of $V$ in $X$.\\
But in general, the blow-up of a coherent ideal in a complex manifold is no longer a smooth complex space.
\item The modification $\tau :\tilde{X} \to X$ associated to the blow-up of a coherent ideal sheaf is always a projective morphism. Then if $X$ is projective (resp. k{\"a}hlerian) so is $\tilde{X}$.\\
\end{enumerate}

When we say that two modifications $\tau_{i} : X_{i}\to X, i = 1,2$ are isomorphic, we mean that there exists an $X-$isomorphism of $X_{1}$ on $X_{2}$.\\

\begin{defn}\label{str fiber product}
Let $X$ be a irreducible  complex space and consider  $\tau_{i} : X_{i}\to X$, for $ i = 1,2$, two  modifications of $X$. We define the {\bf strict fiber product}\footnote{This corresponds to the notion of join in Hironaka's paper : footnote 4 page 193.} of $\tau_{1}$ and $\tau_{2}$ the projection on $X$ of {\bf the} irreducible component of $X_{1}\times_{X}X_{2}$ which surjects on $X$. We shall denote $\tau := \tau_{1}\times_{X,str}\tau_{2} : X_{1}\times_{X,str} X_{2} \to X$ this modification of $X$.
\end{defn}

Note that $X_{1}\times_{X, str}X_{2}$ is irreducible, by definition, and if we have two modifications $f : Y \to X_{1}$ and $g : Y \to X_{2}$ such that $\tau_{1}\circ f = \tau_{2}\circ g := \tau$ then $\tau$ factorizes by a modification $\sigma : Y \to X_{1}\times_{X,str}X_{2}$ induced by $f\times_{X} g$.\\
This notion extends easily to any reduced complex space by taking the union of all irreducible components of the fiber product which surject onto an irreducible component of $X$.

\parag{Example} It is easy to see that $ \tau_{1}\times_{X,str}\tau_{1} : X_{1}\times_{X,str} X_{1} \to X$  is isomorphic to $\tau_{1}$, and that the ``total'' fiber product has an extra component (assuming that $\tau_{1}$ is not injective)  which is the closure of the fiber product of $\tau_{1}^{-1}(C)\times_{C} \tau_{1}^{-1}(C) \setminus \Delta$  where $C$ is the center of $\tau_{1}$ and $\Delta$ the diagonal in $\tau_{1}^{-1}(C)\times_{C} \tau_{1}^{-1}(C)$.$\hfill \square$\\

Remark that, in general, the fiber product of two modifications has extra  irreducible components which project respectively in the  centers of the two modifications.

\begin{prop}\label{finale}
Let $X$ be reduced complex space and let $I$ and $J$ to coherent ideals in $\mathcal{O}_{X}$ such that $V(\mathcal{I})$ and $V(\mathcal{J})$ have no interior point in $X$. Let $\tau_{I} : X_{I} \to X$ and $\tau_{J} : X_{J} \to X$ be the blow-up of $I$ and $J$ in $X$. Let $\tilde{I}$ (resp. $\tilde{J}$) the strict pull-back of $I$ by $\tau_{J}$ (resp. of $J$ by $\tau_{I}$) and let $\sigma_{\tilde{I}}: X_{\tilde{I}} \to X_{J}$ (resp. $\sigma_{\tilde{J}}: X_{\tilde{J}} \to X_{I}$) the blow-up of $\tilde{I}$ (resp. $\tilde{J}$) in $X_{J}$ (resp. $X_{I}$). Then we have :
\begin{enumerate}
\item The modifications $\tilde{J}$ and $\tilde{I}$ are isomorphic.
\item They are isomorphic to the strict fiber product $\tau_{I}\times_{X,str}\tau_{J}$.
\item The blow-up $\tau_{I.J} : X_{I.J} \to X$ is also isomorphic to these modifications, where $I.J$ is the ideal product of $I$ and $J$ in $\mathcal{O}_{X}$.
\end{enumerate}
\end{prop}

The proof of the assertion 3.  will use the following easy lemma.

\begin{lemma}\label{easy}
Let $f : \mathbb{P}_{m}\times \mathbb{P}_{n} \to \mathbb{P}_{m.n+1}$ be the map defined in homogeneous coordinates by
$$(x_{0}, \dots, x_{m}),(y_{0}, \dots, y_{n}) \mapsto (x_{i}.y_{j}), (i,j) \in [0,m]\times [0,n] .$$
This map is a closed embedding.
\end{lemma}

The proof is left to the reader.

\parag{proof of the proposition \ref{finale}} As the strict pull-back of $I$ on $X_{\tilde{I}}$ is locally principal, we have a holomorphic $X-$map $X_{\tilde{I}} \to X_{I}$. But also the pull-back of $\tilde{J}$ by this map is locally principal, so this map lift to a $X-$map to $X_{\tilde{J}}$. This conclude the assertion 1. by symetry.\\
On $X_{\tilde{I}}$ the strict pull-back of $I$ and $J$ are locally principal, so we have a $X-$map to $X_{I}\times_{X,str} X_{J}$. Let  $p : X_{I}\times_{X,str} X_{J} \to X_{J}$ be the projection, the $p_{str}^{*}(\tilde{I})$ is locally trivial on $X_{I}\times_{X,str} X_{J}$, so we have a $X-$map $X_{I}\times_{X,str} X_{J} \to X_{\tilde{I}}$ which is clearly an inverse of the previous one. This proves 2.\\
To prove 3. consider an irreducible small open set $U$ in $X$ on which the ideals $I$ and $J$ are generated respectively by non zero holomorphic functions $f_{0}, \dots, f_{m}$ and $g_{0}, \dots, g_{n}$. Now the modifications $\tau_{I}$ and $\tau_{J}$ are given as follows : take a Zariski dense open set $V \subset U$ on which the maps $F : V \to \mathbb{P}_{m}$ and $G : V \to \mathbb{P}_{n}$ are defined and holomorphic. Then $\tau_{I}^{-1}(U)$ and $\tau_{J}^{-1}(U)$ are respectively the closure of the graphs of $F$ and $G$ in $U \times \mathbb{P}_{m}$ and $U \times \mathbb{P}_{n}$.\\
Now remark that the composition of the map $ \tau_{I}^{-1}(U)\times_{U,str}\tau_{J}^{-1}(U) \to \mathbb{P}_{m}\times \mathbb{P}_{n}$ with the embedding $f$ of the previous lemma gives a proper embedding of $$\tau_{I}^{-1}(U)\times_{U,str}\tau_{J}^{-1}(U) \to \tau_{I.J}^{-1}(U).$$
This gives the point 3. $\hfill \square$\\

An easy corollary of this result is the fact that if the product of two coherent ideals $\mathcal{I}$ and $\mathcal{J}$ is locally principal in $\mathcal{O}_{X}$, then, assuming that $V(\mathcal{I})$ and $V(\mathcal{J})$ have no interior point in $X$,  they are locally principal.

The following example shows that the hypothesis that  $V(\mathcal{I})$ and $V(\mathcal{J})$ have no interior point in $X$ is important.

\parag{Example}  Consider in $\C^{4}$ with coordinates $(x, y, u, v)$ the ideal generated by $$ x.u - y.v, \quad x.v - y.u, \quad y.u - y.v .$$
 It define a surface $S$ which is the union of three $2-$planes :
$$ P_{1}:= \{  u = v = 0\}, \quad P_{2} := \{ x = y = 0\}, \quad P_{3} := \{ x = y \quad  { \rm and} \quad  u = v \} .$$
The intersection $P_{1}\cap P_{2}$ is reduced to $\{0\}$ but $P_{3}$ meets $P_{1}$ and $P_{2}$ respectively in the lines $\{ u = v = x - y = 0\}$ and $\{u - v = x = y = 0 \}$.
Moreover, it is easy to verify that the given ideal is reduced, so $S$ is a reduced complex surface (but not irreducible). It is immediate to see that the ideals $\mathcal{I} := (x, y)$ and $\mathcal{J} := (u, v)$ are not principal in the ring $\mathcal{O}_{S,0} $ but that their product is principal and generated by the element $\varphi := x.u = x.v = y.u = y.v $ in the ring $\mathcal{O}_{S,0} $.\\

\parag{Remark} If we consider a complex manifold $X$ and coherent ideals $I_{1}, \dots, I_{k}$ in $\mathcal{O}_{X}$ such that the $V(\mathcal{I}_{i}), i \in [1,k]$ has no interior point in $X$,  to compute the blow-up of the ideal
$$ J := I_{1}^{p_{1}}.\dots.I_{k}^{p_{k}} \quad {\rm where} \quad p_{1}, \dots, p_{k} \quad {\rm are \ positive \ integers} $$
we can simply blow-up successively $I_{1}, \dots, I_{k}$ ignoring the exponents $p_{j}, j \in [1,k]$.\\

We shall also use the following lemma later on.

\begin{lemma}\label{eclat; induit}
Let $X$ be a reduced complex space and $Z \subset X$ be a closed analytic subset. Let $I $ a coherent in $\mathcal{O}_{X}$ and assume that $V(I)$, the locus of common zeros of sections of $I$, has no interior point in $X$ and in  $Z$. Then we have a commutative diagram of blow-up
$$\xymatrix{Z_{I} \ar[r]^{i_{I}} \ar[d]_{\sigma_{I}}& X_{I} \ar[d]_{\tau_{I}} \\ Z \ar[r]^{i} & X } $$
where $i_{I}$ is a proper embedding, and where $\sigma_{I} : Z_{I} \to Z$ is the blow-up in $Z$ of the image in $\mathcal{O}_{Z}$ of the ideal $I$.
\end{lemma}

\parag{proof} Let $U$ be an irreducible open set in $X$ such that $Z \cap U$ is also irreducible and such that $I$ is generated on $U$ by holomorphic functions $f_{0}, \dots, f_{m}$. Let $F : U\setminus V(I) \to \mathbb{P}_{m}$ be the map defined by $f_{0}, \dots, f_{m}$. Then we have
 $$\tau_{I}^{-1}(U) \simeq \overline{graphe \  F} \subset U \times \mathbb{P}_{m}, \  {\rm and} \quad \sigma_{I}^{-1}(U \cap X) \simeq \overline{graphe \ F_{\vert X}} \subset (U\cap X) \times \mathbb{P}_{m}.$$
  This is enough to conclude.$\hfill \square$\\
  
  \section{Smoothness over $Q(u)$}

Now we come back to the example of [H.61]. An important point is to blow-up ideals which are product of (simple) ideals in order to use the proposition \ref{finale}. The ideal which appears in the next statement will be crucial for the construction of the manifold $V_{0}$.

\begin{lemma}\label{H1}{\rm [cf. [H.61]]}
Let $(x_{1}, x_{2},x_{3})$ the coordinates in $\C^{3}$. Then we have the equality of ideals in $\mathcal{O}_{\C^{3}}$ 
\begin{align*}
& \mathcal{H} :=  (x_{1}.x_{2}, x_{2}.x_{3}, x_{1}.x_{3}).(x_{1}, x_{2}.x_{3}).(x_{2}, x_{1}.x_{3})^{2}.(x_{3}, x_{1}.x_{2})^{2} \\
& \quad =  (x_{2}, x_{3})^{5}\cap (x_{1},x_{3})^{4}\cap (x_{1}, x_{2})^{4}\cap (x_{1}, x_{2}, x_{3})^{7}.
\end{align*}
\end{lemma}

\parag{Proof} It is enough to compare monomials in these ideals. A necessary and sufficient condition for $x_{1}^{a}.x_{2}^{b}.x_{3}^{c}$ to be in the right-handside is given by the following inequalities :
\begin{equation*} a+ b+ c \geq 7, \quad a+b \geq 4 , \quad a+ c \geq 4 , \quad {\rm and} \quad b+c \geq 5 . \tag{@}
\end{equation*}
So listing the cases for $a = 0$ to $a = 4$ gives the following generators for the right-handside :
\begin{align*}
& m_{0} := x_{2}^{4}.x_{3}^{4}, \quad  m_{1} := x_{1}.x_{2}^{3}.x_{3}^{3}, \quad m_{2} := x_{1}^{2}.x_{2}^{3}.x_{3}^{2},  \quad m_{2}' := x_{1}^{2}.x_{2}^{2}.x_{3}^{3}, \\
& m_{3} = x_{1}^{3}.x_{2}^{4}.x_{3}, \quad  m_{3}' := x_{1}^{3}.x_{2}.x_{3}^{4}, \quad m_{4} := x_{1}^{4}.x_{2}^{5}, \quad m_{4}' := x_{1}^{4}.x_{3}^{5}.
\end{align*} 
Now it is rather easy to verify (and this will be indicate in the verification of the opposite inclusion) that each of these monomials is in the right-handside:
\begin{align*}
& m_{0} = x_{2}^{4}.x_{3}^{4} = (x_{2}.x_{3}).(x_{2}.x_{3}).x_{2}^{2}.x_{3}^{2} \\
& m_{1} =  x_{1}.x_{2}^{3}.x_{3}^{3} = (x_{2}.x_{3}).x_{1}.x_{2}^{2}.x_{3}^{2} \\
& m_{2} = x_{1}^{2}.x_{2}^{3}.x_{3}^{2} = (x_{1}.x_{2}).x_{1}.x_{2}^{2}.x_{3}^{2} \\
& m_{2}' = x_{1}^{2}.x_{2}^{2}.x_{3}^{3} = (x_{1}.x_{3}).x_{1}.x_{2}^{2}.x_{3}^{2} \\
& m_{3} = x_{1}^{3}.x_{2}^{4}.x_{3} = (x_{1}.x_{2}).x_{1}.x_{2}^{2}.(x_{1}.x_{2}.x_{3}) \\
& m_{3}' = x_{1}^{3}.x_{2}.x_{3}^{4} = (x_{1}.x_{3}).x_{1}.(x_{1}.x_{2}.x_{3}).x_{3}^{2} \\
& m_{4} =  x_{1}^{4}.x_{2}^{5} = (x_{1}.x_{2}).x_{1}.x_{2}^{2}.x_{3}^{2} \\
& m_{4}' =  x_{1}^{4}.x_{3}^{5} =  (x_{1}.x_{3}).x_{1}.(x_{1}.x_{3})^{2}.x_{3}^{2} \\
\end{align*}

For $a \geq 5$, as we have $b + c \geq 5$, it is easy to check that we obtain multiples of the eight monomials above : if $b$ or $c$ is $0$, a multiple of $m_{4}$ or $m_{4}'$ respectively, if $b$ or $c$ is equal to $1$ a multiple of $m_{3}$ or $m_{3}'$ respectively  and if  $b$ or $c$ is at least $2$ a multiple of $m_{2}$ or $m_{2}'$.\\

It is a little more painful to verify the opposite inclusion because there are {\it a priori} 54 monomials in the generator of the left-handside. \\
We shall use the symbol $(@)$ to indicate when we find one of the eight monomials above in the list of these 54 monomials.\\

The maximal degree in $x_{1}$ for such a monomial is $6$. In degree $6$ there are only two: 
$$(x_{1}.x_{2}).x_{1}.(x_{1}.x_{3})^{2}.(x_{1}.x_{3})^{2}  = x_{1}^{6}.x_{2}^{3}.x_{3}^{2}$$
 and the one obtain by exchanging $x_{2}$ and $x_{3}$ (we note this by (ex-2-3)). They are multiples of $m_{2}$ and $m_{2}'$ respectively.\\
 
 In degree $5$ in $x_{1}$ we have
 \begin{align*}
 & (x_{2}.x_{3}).x_{1}.(x_{1}.x_{3})^{2}. (x_{1}.x_{2})^{2} =  x_{1}^{5}.x_{2}^{3}.x_{3}^{3} \quad \in (m_{1}) \\
 & (x_{1}.x_{2}).(x_{2}.x_{3}).(x_{1}.x_{3})^{2}. (x_{1}.x_{2})^{2} =  x_{1}^{5}.x_{2}^{4}.x_{3}^{3} \quad {\rm (ex-2-3)} \quad  \in (m_{1}) \\
 &  (x_{1}.x_{2}).x_{1}.(x_{1}.x_{2}.x_{3}).(x_{1}.x_{2})^{2} =  x_{1}^{5}.x_{2}^{4}.x_{3} \quad {\rm (ex-2-3)} \quad  \in (m_{3}) \quad {\rm and} \quad (m_{3}') \\
 &  (x_{1}.x_{2}).x_{1}.(x_{1}.x_{3})^{2}.(x_{1}.x_{2}.x_{3}) = x_{1}^{5}.x_{2}^{2}.x_{3}^{3} \quad {\rm (ex-2-3)} \quad  \in (m_{2}') \quad {\rm and} \quad (m_{2}) \\
  \end{align*}
  
  Now consider the monomial of degree $4$ in $x_{1}$. First consider these where we have two $x_{1}^{2}$ coming from one of the four terms.
  
  \begin{align*}
  & (x_{2}.x_{3}).(x_{2}.x_{3}).(x_{1}.x_{3})^{2}. (x_{1}.x_{2})^{2} = x_{1}^{4}.x_{2}^{4}.x_{3}^{4} \quad \in (m_{1}) 
  & 
  \end{align*}
   
   Now if we have only one $x_{1}^{2}$ coming for one four terms 
   
  \begin{align*}
  &  (x_{2}.x_{3}).x_{1}.(x_{1}.x_{2}.x_{3}).(x_{1}.x_{2})^{2} =  x_{1}^{4}.x_{2}^{4}.x_{3}^{2}  \quad  \in (m_{2}) \\
  &  (x_{1}.x_{2}).(x_{2}.x_{3}).(x_{1}.x_{2}.x_{3}).(x_{1}.x_{2})^{2} =  x_{1}^{4}.x_{2}^{5}.x_{3}^{2} \quad {\rm (ex-2-3)} \quad  \in (m_{2})  \quad {\rm and} \quad (m_{2}')\\
  &  (x_{1}.x_{2}).x_{1}.x_{2}^{2}.(x_{1}.x_{2})^{2} =  x_{1}^{4}.x_{2}^{5}  \quad  {\rm (ex-2-3)} \quad  \in (m_{4}) \quad {\rm and} \quad (m_{4}') \quad  (@) + (@)\\ 
  &   (x_{1}.x_{2}).x_{1}.(x_{1}.x_{3})^{2}.x_{3}^{2} =  x_{1}^{4}.x_{2}.x_{3}^{4} \quad  {\rm (ex-2-3)} \quad  \in (m_{3}')  \quad {\rm and} \quad (m_{3}) \\  
  &   (x_{2}.x_{3}).x_{1}. (x_{1}.x_{3})^{2}.(x_{1}.x_{2}.x_{3}) =  x_{1}^{4}.x_{2}^{2}.x_{3}^{4} \quad  \in (m_{2}')\\
  &    (x_{1}.x_{2}). (x_{2}.x_{3}).(x_{1}.x_{3})^{2}.(x_{1}.x_{2}.x_{3}) =  x_{1}^{4}.x_{2}^{3}.x_{3}^{4} \quad {\rm (ex-2-3)}  \quad  \in (m_{1}) \quad {\rm and} \quad  (m_{1})\\
  \end{align*}
   
 If we choose one $x_{1}$ in each term we find only
\begin{align*}
& (x_{1}.x_{2}).x_{1}.(x_{1}.x_{2}.x_{3}).(x_{1}.x_{2}.x_{3}) = x_{1}^{4}.x_{2}^{3}.x_{3}^{2} \quad {\rm (ex-2-3)} \quad  \in (m_{2}) \quad {\rm and} \quad (m_{2}')
\end{align*}

In degree $3$ in $x_{1}$ if we take $x_{1}^{2}$ in one term we will get
\begin{align*}
&  (x_{1}.x_{2}).(x_{2}.x_{3}).x_{2}^{2}.(x_{1}.x_{2})^{2} =  x_{1}^{3}.x_{2}^{6}.x_{3} \quad {\rm (ex-2-3)} \quad  \in (m_{3}) \quad {\rm and} \quad  (m_{3}') \\
&   (x_{2}.x_{3}).x_{1}.x_{2}^{2}.(x_{1}.x_{2})^{2} =  x_{1}^{3}.x_{2}^{5}.x_{3} \quad  \in (m_{3})  \\
&    (x_{2}.x_{3}). (x_{2}.x_{3}).(x_{1}.x_{2}.x_{3}).(x_{1}.x_{2})^{2} =  x_{1}^{3}.x_{2}^{5}.x_{3}^{3} \quad \in (m_{1}) \\
&  (x_{1}.x_{2}).(x_{2}.x_{3}).(x_{1}.x_{3})^{2}.x_{3}^{2} =  x_{1}^{3}.x_{2}^{2}.x_{3}^{5} \quad {\rm (ex-2-3)} \quad \in  (m_{2}' )\quad {\rm and} \quad (m_{2}) \\
&  (x_{2}.x_{3}).x_{1}.(x_{1}.x_{3})^{2}.x_{3}^{2} =  x_{1}^{3}.x_{2}.x_{3}^{5}  \quad  \in  (m_{3}') \\
&   (x_{2}.x_{3}). (x_{2}.x_{3}).(x_{1}.x_{3})^{2}.(x_{1}.x_{2}.x_{3}) = x_{1}^{3}.x_{2}^{3}.x_{3}^{5}  \quad  \in  (m_{1})
\end{align*}

For the degree $3$ in $x_{1}$, if we take at most one $x_{1}$ in each term, we have four places where to avoid $x_{1}$ :
\begin{align*}
& (x_{2}.x_{3}).x_{1}.(x_{1}.x_{2}.x_{3}).(x_{1}.x_{2}.x_{3}) = x_{1}^{3}.x_{2}^{3}.x_{3}^{3} \quad  \in (m_{1}) \\
& (x_{1}.x_{2}).(x_{2}.x_{3}).(x_{1}.x_{2}.x_{3}).(x_{1}.x_{2}.x_{3}) = x_{1}^{3}.x_{2}^{4}.x_{3}^{3} \quad {\rm (ex-2-3)}  \quad  \in (m_{1}) \\
& (x_{1}.x_{2}).x_{1}.x_{2}^{2}.(x_{1}.x_{2}.x_{3}) = x_{1}^{3}.x_{2}^{4}.x_{3} \quad {\rm (ex-2-3)}  \quad  \in (m_{3} )\quad {\rm and} \quad  (m_{3}') \quad (@) + (@)\\
& (x_{1}.x_{2}).x_{1}.(x_{1}.x_{2}.x_{3}).x_{3}^{2} = x_{1}^{3}.x_{2}^{2}.x_{3}^{3} \quad {\rm (ex-2-3)}  \quad  \in (m_{2}') \quad {\rm and} \quad  (m_{2})\\
\end{align*}
In degree $2$ there are $8$ choices, modulo exchanging $x_{2}$ and $x_{3}$:
\begin{align*}
&  (x_{1}.x_{2}).x_{1}.x_{2}^{2}.x_{3}^{2} = x_{1}^{2}.x_{2}^{3}.x_{3}^{2}\quad {\rm (ex-2-3)} \quad  \in (m_{2}) \quad {\rm and} \quad (m_{2}') \quad (@) + (@)  \\
&  (x_{1}.x_{2}).(x_{2}.x_{3}).(x_{1}.x_{2}.x_{3}).x_{3}^{2} = x_{1}^{2}.x_{2}^{3}.x_{3}^{4}\quad {\rm (ex-2-3)}  \quad \in (m_{1}) \\
&   (x_{1}.x_{2}).(x_{2}.x_{3}).x_{2}^{2}.(x_{1}.x_{2}.x_{3}) = x_{1}^{2}.x_{2}^{5}.x_{3}^{2}\quad {\rm (ex-2-3)} \quad  \in (m_{2}) \quad {\rm and} \quad  (m_{2}') \\ 
&  (x_{2}.x_{3}).x_{1}.(x_{1}.x_{2}.x_{3}).x_{3}^{2} = x_{1}^{2}.x_{2}^{2}.x_{3}^{4} \quad  \in (m_{2}')  \\
&   (x_{2}.x_{3}).x_{1}.x_{2}^{2}.(x_{1}.x_{2}.x_{3}) = x_{1}^{2}.x_{2}^{4}.x_{3}^{2} \quad  \in  (m_{2}) \\
&   (x_{2}.x_{3}).(x_{2}.x_{3}).(x_{1}.x_{2}.x_{3}).(x_{1}.x_{2}.x_{3}) = x_{1}^{2}.x_{2}^{4}.x_{3}^{4} \quad  \in (m_{1}) \\
&    (x_{2}.x_{3}).(x_{2}.x_{3}).(x_{1}.x_{3})^{2}.x_{3}^{2} =  x_{1}^{2}.x_{2}^{2}.x_{3}^{6} \quad  \in (m_{2}') \\ 
&    (x_{2}.x_{3}).(x_{2}.x_{3}).x_{2}^{2}.(x_{1}.x_{2})^{2} = x_{1}^{2}.x_{2}^{6}.x_{3}^{2} \quad  \in (m_{2})  \\
\end{align*}
In degree $1$ there are $4$ choices, modulo exchanging $x_{2}$ and $x_{3}$:
\begin{align*}
&   (x_{1}.x_{2}). (x_{2}.x_{3}).x_{2}^{2}.x_{3}^{2} =  x_{1}.x_{2}^{4}.x_{3}^{3}\quad {\rm (ex-2-3)} \quad  \in  (m_{1}) \\
&   (x_{2}.x_{3}).x_{1}.x_{2}^{2}.x_{3}^{2} =  x_{1}.x_{2}^{3}.x_{3}^{3} \quad  \in (m_{1}) \quad (@) \\
&    (x_{2}.x_{3}).(x_{2}.x_{3}).(x_{1}.x_{2}.x_{3}).x_{3}^{2} =  x_{1}.x_{2}^{3}.x_{3}^{5} \quad  \in (m_{1}) \\
&     (x_{2}.x_{3}).(x_{2}.x_{3}).x_{2}^{2}.(x_{1}.x_{2}.x_{3}) =  x_{1}.x_{2}^{5}.x_{3}^{3} \quad  \in  (m_{1}) \\
\end{align*}
For the degree $0$ there is only one possibility : $x_{2}^{4}.x_{3}^{4}$ which is $m_{0} \ (@)$.
So the verification for the $54$ monomial is over.$\hfill \square$\\
\bigskip

Let $I_{C_{i}}$ be the reduced ideal of the smooth curve $C_{i}$ for $i = 1,2,3$. Then 
we have the following consequence of the Lemma \ref{H1}.

\begin{cor}
The ideal $\mathcal{H} $ defined  in \ref{H1} admits the following alternative description in the complement of the origin of 
$\C^{3}$
\begin{enumerate}
\item[\rm (1)] if $x_{1} \not= 0$, $\mathcal{H} $ is equal to $(x_{2}, x_{3})^{5} = I_{C_{1}}^{5} $,
\item[\rm (2)] if $ x_{2} \not= 0$, $\mathcal{H} $ is equal to $(x_{1}, x_{3})^{4} = I_{C_{2}}^{4}$,
\item[\rm (3)] if  $x_{3}\not= 0 $, $\mathcal{H} $ is equal to $(x_{1}, x_{2})^{4} = I_{C_{3}}^{4}$.
\end{enumerate}
\end{cor}

We blow-up the product  ideal $I_{C_{1}}^{5}.I_{C_{2}}^{4}.I_{C_{3}}^{4}$ in 
$\mathbb{P}_{3} \setminus \{P\}$. In a neighborhood of $P$ we glue 
the blow-up of $\C^{3}$ 
along the ideal $\mathcal{H} $ in \ref{H1}; the resulting reduced complex space will be our central fiber $V_{0}$. It is possible to perform such a surgery 
precisely because the points (1)--(3) of the preceding corollary and the remark before the lemma \ref{eclat; induit}.\\
Note that the remark 2 before the definition \ref{str fiber product} gives that $V_{0}\setminus \{P, Q\}$  is smooth.\\ 

Formally, the space $V_0$ is equal to the fibered product of the two blow-ups near $P$, and it coincides with the 
blow-up of $\bP_3$ ``far" from this point. 
The smoothness of $V_{0}$ is reduced to check the smoothness of the blow-up of $\C^{3}$ along the ideal $\mathcal{H}$ of Lemma \ref{H1} near $0 \ (= P)$ and $Q$ which belongs to the chart $U_{1}$.\\

Note also that when Hironaka considers the curve $C_{3}(u)$ for $u \not= 0$, the 
ideal $\mathcal{H}$  is globally defined on $\mathbb{P}_{3}$ using the corollary above. Near $Q$ it coincides with the intersection $I_{C_{1}}^{5}\cap I_{C_{2}}^{4}$ and near $Q'(u) \not= Q$, the second intersection point of $C_{1}\cap C_{3}(u)$, to $I_{C_{1}}^{5}\cap I_{C_{3}(u)}^{4}$. \\
This gives the projectivity of the $u-$fiber of the family.\\

 In order to get the smoothness at $Q$ or at $Q'(u)$ in this case we shall use the next lemma.

\begin{lemma}\label{H1bis}
Let $(x_{1}, x_{2}, x_{3})$ the standard coordinates in $\C^{3}$. Then we have the equality of ideals :
$$ (x_{2}, x_{3})^{5} \cap (x_{1}, x_{3})^{4} = (x_{1}.x_{2}, x_{3})^{4}.(x_{2}, x_{3}) .$$
\end{lemma}

\parag{proof} The generators of the right handside are the monomials $x_{1}^{a}.x_{2}^{a+1}.x_{3}^{b}$ and $x_{1}^{a}.x_{2}^{a}.x_{3}^{b+1}$ where $a + b = 4$. They are cleraly in the left handside. \\
Conversely, the left handside is generates by monmials $x_{1}^{\alpha}.x_{2}^{\beta}.x_{3}^{\gamma}$ where we have $\alpha + \gamma = 4$ and $\beta + \gamma = 5$. To show that there are in the right handside, it is enough to show that for any such $(\alpha, \beta, \gamma)$ we can find non negative integers $(a,b)$ such that $a + b \geq 4$ with
\begin{enumerate}
\item  either $a \leq \alpha, \quad a+1 \leq \beta$ and $b \leq \gamma$, \hfill  case 1)
\item or  $a \leq \alpha, \quad a+1 \leq \beta$ and $b \leq \gamma$, \hfill  case 2)
\end{enumerate}
This is obtained as follows:
\begin{align*}
&  \beta = 0 \quad \gamma = 5 \quad \alpha = 0 \quad   then \quad  a = 0  \quad  b = 4  \quad  case  \ 2) \\
& \beta = 1 \quad  \gamma = 4 \quad \alpha = 0 \quad    then \quad a = 0  \quad  b = 4  \quad case \  1) \\
& \beta = 2 \quad  \gamma = 3 \quad \alpha = 1 \quad    then \quad a = 1  \quad  b = 3  \quad case \ 1) \\
& \beta = 3 \quad  \gamma = 2 \quad \alpha = 2 \quad    then \quad a = 2  \quad  b = 2  \quad case\ 1) \\
& \beta = 4 \quad  \gamma = 1 \quad \alpha = 3 \quad    then \quad a = 3  \quad  b = 1  \quad case \ 1) \\
& \beta = 5 \quad  \gamma = 0 \quad \alpha = 4 \quad    then \quad a = 4  \quad  b = 0  \quad case \ 1) 
\end{align*}
This completes the proof.$\hfill \square$\\

The smoothness of the blow-up at $Q$ and $Q'(u)$ for $u \not= 0$ is then an easy corollary of this lemma because after the blow-up of the ideal $(x_{2},x_{3})$ which corresponds to $C_{1}$, the the strict pull-back of the ideal is either principal or is the ideal of a smooth curve. Hence the blow-up of the ideal of the lemma is smooth.

\section{Smoothness over $P$}

We shall consider the blow-up  $X$  of  the ideal $(x_{1}.x_{2}, x_{2}.x_{3}, x_{3}.x_{1})$ in $\C^{3}$ and then the blow up $Y$ of the ideal $(x_{3}, x_{1}.x_{2} ).(x_{1}, x_{2}.x_{3} ).(x_{2}, x_{3}.x_{1} )$ in $X$. This corresponds to the blow-up of the ideal described in the lemma \ref{H1} which gives the first piece in the construction of the manifold $V_{0}$.\\

In $\C^{3}\times \mathbb{P}_{2}$ we look for the irreducible component of  the analytic subset
$$ Z := \{((x_{1}, x_{2}, x_{3}), (a,b,c)) \ / \  x_{1}. x_{2}.b =  x_{2}. x_{3}.a \quad x_{1}.x_{2}.c = x_{3}.x_{1}.a \quad x_{2}.x_{3}.c = x_{3}.x_{1}.b \} $$
which dominates the graph of the corresponding map. By circular permutation on $x_{1}, x_{2}, x_{3}$ we can restrict our study to the chart $a \not= 0$ in $\mathbb{P}_{2}$. Then let $u: = b/a$ and $v: = c/a$ be the corresponding coordinates. As $x_{1}, x_{2}, x_{3}$ are not zero at the generic point of the graph, the equations simplify in
$$ x_{3} = u.x_{1} \quad x_{3} = v.x_{2} \quad x_{2}.v = x_{1}.u  $$
and this show that this Zariski open  set  $X'$ in $X$ is isomorphic to the hypersurface $H: = \{ x_{2}.v = x_{1}.u  \} \subset \C^{4}$ where the coordinates in $\C^{4}$ are $(x_{1}, x_{2}, u, v)$.
So we have an unique singular (Morse) point: the origin. \\

Note that the pre-image of  $0 $ in $X'$ is the plane defined by $x_{1} = x_{2} = 0$.\\
In this chart, isomorphic to the hypersuface $H$  note 
 $$P_{1}: = \{ x_{2} = u = 0\}, \quad  P_{2} := \{ x_{1} = v = 0\},\quad  P_{3} := \{x_{1} = x_{2} = 0\} .$$
  We have $P_{1}\cap P_{2} = \{0\}, P_{1}\cap P_{3} = \{x_{1} = x_{2} = u = 0 \}, P_{2} \cap P_{3} = \{ x_{1} = x_{2} = v = 0\}$. The pull back of the curves $C_{1} = \{x_{1} = x_{2} = 0\}$ is $P_{1}\cup P_{3}$,  the pull-back  of the curve $C_{2} = \{ x_{1} = x_{3} = 0 \}$ is $P_{2}\cup P_{3}$ and  the pull-back  of the curve $C_{3} = \{x_{1} = x_{2} = 0 \}$ is $P_{3}$. See the picture $0$.\\

Remark now that the ideals $(x_{1}, x_{2}.x_{3})$ and $(x_{2}, x_{1}.x_{3})$ are principal on $X'$ as we have $x_{3} = u.x_{1} = v.x_{2}$. So the only blow-up to perform in $X'$ is now the blow up of  the ideal $(x_{3}, x_{1}.x_{2})$. Using $x_{3} = u.x_{1}$ this is the same as the blow-up of $(u, x_{2})$. Then it is given in $X'\times \mathbb{P}_{1}$ by the equation $u.\beta = x_{2}.\alpha$. Let $z := \beta/\alpha$ for the first chart in $\mathbb{P}_{1}$. Then we have in this chart  $X_{1}''$ of this blow-up
$$ x_{2} = u.z, \quad  x_{2}.v = x_{1}.u \quad x_{3} = v.x_{2} = u.x_{1} $$
which give, as $x_{2}$ is generically non zero,  $x_{1} = v.z , \quad x_{2} = u.z ,\quad  x_{3} = u.v.z $ \ and then we have a copy of $\C^{3}$ with coordinates $(u, v, z)$. The computations on the other charts are analogous as $x_{1}$ and $x_{2}$ have symetric roles in $X'$ and in the ideal blown up.\\

So we conclude that the blow-up of $\C^{3}$ by the ideal
 $$(x_{1}.x_{2}, x_{2}.x_{3}, x_{3}.x_{1}).(x_{1}, x_{2}.x_{3}).(x_{2}, x_{3}.x_{1} )^{2}.(x_{3}, x_{1}.x_{2} )^{2}$$
  is a smooth quasi-projective manifold.\\

Remark that the pull-back of the maximal ideal of $\C^{3}$  on $X'$ is the ideal $(x_{1}, x_{2})$ in $X'$ which is a copy of $\C^{2}$ with coordinnates $(u, v)$. Now, it becomes the ideal $z.(u,v)$ in $X''$ which is the union of a $\C^{2} := \{ z = 0 \}$ and a transversal line which has only its intersection with this plane over the origin in $\C^{3}$. \\
So our exceptional divisor is just this plane $\C^{2}$. Now, as $u = b/a$ and $v := c/a$ it is easy to see that the global exceptional divisor for the final modification of $\C^{3}$  is a copy of $\mathbb{P}_{2}$. This divisor is called  $E_{0}$ in the pictures 1.\\
The intersection, in this chart, of this $\mathbb{P}_{2}$ which  is view as  $\C^{2} = \{ z = 0 \}$, with the strict transform of the curve $\{x_{1}= x_{3} = 0\}$ is  $\{z = v = 0\}$, and with the strict transform  the curve $\{x_{2} = x_{3} = 0\}$ is $\{ z = u = 0 \}$. \\
Using the six charts for the blow-up, we find that each strict transform of the curve $\{x_{3} = x_{1} = 0 \}$  cuts (generically transversally) the exceptional divisor in $3$ lines in general position. This corresponds to the picture 1. The visible part in the chart $X''$ is given by the picture 2. \\

Near the point $Q$ we simply blow-up $C_{1}$ and then the strict transform of $C_{2}$ and $C_{3}$. This is equivalent to blow-up the product of the reduced ideals. As $C_{1}$ is smooth as well as the strict transforms of $C_{2}$ and $C_{3}$ which are disjoint  in $U'_{1}$, the final blow-up is smooth and patches with the blow-up near $P$.\\
The pictures 3 explains the situation near the point $Q$ and the picture 4 the global situation in $\mathbb{P}_{3}$ before the blowing-up.

\parag{Why $V_{0}$ is not k{\"a}hlerian} We can read on pictures 1 and 3 the following algebraic equivalence of curves in $V_{0}$ :
\begin{align*}
& L_{1} \sim L'_{1} + L_{2} + L_{3} \quad\quad {\rm on \ picture \ 3} \\
& L_{1} \sim L_{1,2} + L_{0,1} + L_{1,3} \quad {\rm on \ picture \ 1} \\
& L_{2} \sim L_{2,3} + L_{1,2} + L_{0,2} \quad {\rm on \ picture \ 1} \\
& L_{3} \sim L_{3,1} + L_{2,3} + L_{0,3} \quad {\rm on \ picture \ 1} 
\end{align*}
and this implies,
$$ L_{0,1} \sim L'_{1} + 2L_{2,3} +  L_{0,2} + L_{0,3} .$$
Now the curve $L_{0,1}, L_{0,2}, L_{0,3}$ are lines in $E_{2} \simeq \mathbb{P}_{2}$ so they are algebraically equivalent. This implies that $L'_{1} + 2L_{2,3} + L_{0,2} \sim 0$ in $V_{0}$ and prove our claim.\\

Remark that for $u \not=  0$ the curve $C_{3}(u)$ no longer meets $C_{1}$ at $Q$ so the first relation above becomes $L_{1} \sim L'_{1} + L_{2}$ and now the computation analog to the computation above only gives $L'_{1} + L_{2,3} \sim L_{1,3}$ which does not contradict the projectivity of the fiber at $u$ of the family. The projectivity of the $u-$fibers when $u \not= 0$ is consequence of the fact that we obtain the blow-up of a global coherent ideal in $\mathbb{P}_{3}$. Moreover this shows that the restriction of the projection of the  family over the open set $\{u \not= 0\}$ is a locally projective morphism.

\section{Simpler examples}

 Let $\mathcal{B}$ a ball in $\C^{3}$ with center $0$ and let $C_{1}$ and $C_{2}$ be smooth connected curves defined in $\mathcal{B}$ and meeting transversaly at $0$ and nowhere else. Let $X_{1}$ the blow-up of $\mathcal{B}$ along the reduced ideal of $C_{1}$ and denote by $C'_{2}$ the strict transform of $C_{2}$ in $X_{1}$ and $E_{1}$ the exceptional divisor in $X_{1}$. If $\tau_{1} : X_{1} \to \mathcal{B}$ is the blow-up map, we have $E_{1} = \tau_{1}^{-1}(C_{1})$ and it is a reduced smooth divisor isomorphic to the projectivized normal bundle to $C_{1}$ in $\mathcal{B}$. The curve $C'_{2}$ meets $E_{1}$ transversally in a (unique) point $\tilde{P}$ which belongs to the curve $ A := \tau_{1}^{-1}(0)$ which is a smooth rational curve (the projective space of the fiber at $0$ of the normal bundle of $C_{2}$ at the origin).\\
Now we blow-up $C'_{2}$ in $X'$ to obtain a modification $\sigma_{1,2} : X_{1,2} \to X_{1} $. Then the strict transform of $C'_{2}$ is a smooth divisor $E_{2}$ in $X_{1,2}$ and the strict transform $E'_{1}$ of $E_{1}$ is, via the map $\sigma_{1,2}$, isomorphic to the blow-up of $E_{1}$ at the point $\tilde{P}$. Let $B$ be the smooth rational curve $\sigma_{1,2}^{-1}(\tilde{P})$. Denote by $Z_{1}$ the generic fiber on $E_{1}$  (or of $E'_{1}$) over $C_{1}$. Then in $X_{1,2}$ we have $Z_{1} \sim A' + B$, where $A'$ is the strict transform of $A$ in $E'_{1}$  by the blow-up of $\tilde{P}$ in $E_{1}$. So $A'\cup B$ is the pull back of $A$ in $X_{1,2}$. Remark also that $B$ is homologous in $X_{1,2}$ to the generic fiber $Z_{2}$ of $E_{2}$ over $C'_{2}$ (or $C_{2}$).\\

Let $C$ be a connected curve in $\mathbb{P}_{3}$ with an ordinary double point $P$ and which is smooth outside $P$. Consider now $\mathcal{B}$ as an open neighbourhood of $P$ in $\mathbb{P}_{3}$ and assume that $C_{1}$ and $C_{2}$ are the two branches of the curve $C$ near $P$. In the complex manifold $ V := \mathbb{P}_{3}\setminus \overline{\mathcal{B}'}$, where $\mathcal{B}' \subset \subset \mathcal{B}$ is an open ball with a smaller radius, consider the blow-up $\sigma :Y_{\infty} \to V$ of the smooth connected curve $C \cap V$. Then on the open set $\mathcal{B}\setminus \overline{\mathcal{B}'}$ we have a natural identification between $\sigma$ and $\tau_{1,2}$, because, outside the origin, it is the same thing to blow-up successively $C_{1}$ and $C_{2}$ or to blow-up $C$.\\
So we can glue these two maps to obtain a complex manifold $Y$ and a modification $\tau : Y \to \mathbb{P}_{3}$. Now in $Y$ we have $Z_{1} \sim A' + B$ and also $Z_{2} \sim B$. But the generic fiber of the gluing of $E'_{1} \cup E_{2}$ with the blow-up of $C$ in $V$ implies that $Z_{1} = Z_{2}$. Then the curve $A'$ is homologuous to $0$ in $Y$. See picture 5.\\

If we consider a holomorphic family $(C_{s})_{s \in D}$ of curves parametrized by the unit disc in $\C$ such that $C_{0} = C$ and such that for $s \not= 0$ the curve $C_{s}$ has no singular point, it seems interesting to blow-up (in a suitable way) the graph $\Gamma$ of this family of curve in $D \times \mathbb{P}_{3}$ in order to obtain an holomorphic family of complex manifolds $(Y_{s})_{s \in D}$ such that for $s \not= 0$ the manifold $Y_{s}$ is projective, and such that $Y_{0}$ is not K{\"a}hler. Of course this would give a much more simple example of such a family than the previous example of H. Hironaka from his paper [H.61] discussed above.\\

We explain next the reason why this does not work.\\
The graph $\Gamma \subset D \times \mathbb{P}_{3}$ is a complex submanifold excepted at the point $\{0\}\times \{P\}$ where we have a normal crossing point for a surface in $\C^{4}$. In order to perform a blowing-up ``in family'' for the curve $C_{s}$, we want to separate the two branches of this surface near $\{0\}\times \{P\}$. Of course, this is possible set theoretically, but we want that the fiber at \ $s = 0$ \ of this ''family blow-up'' will give us the previous construction at $t=0$. This essentially means that we want that the projection $\Gamma \to D$, which is is flat  when we consider the reduced structure for $\Gamma$, has a {\bf reduced} fiber at $s = 0$. But we shall show now that  { \bf this is not the case} in our situation. So it is not possible to separate the two irreducible branches of the fiber at $0$  near $0$ keeping the non reduced structure of this fiber !

\begin{lemma}\label{hyp. flat}
Let $X \subset D \times M$ be an analytic subset where $D \subset \C$ is the unit disc and $M$ a complex manifold. Assume that $X$ is reduced and flat on $D$. Let $X_{0}$ be the fiber at $0$ of the projection $\pi : X \to D$ with its ``fiber structure'' meaning that $\I_{X_{0}} := \I_{X}\big/s.\I_{X} \subset \mathcal{O}_{M}$. Assume that this fiber is contained in a smooth hypersurface $H_{0}$ in $M$ near a point $ p \in X_{0}$\footnote{This precisely mean that there exists a holomorphic function $f_{0}$ in an open neighbourhood $U$ of $p$ such that $(df_{0})_{p} \not= 0$ and $f_{0}$ is a section of $\I_{X_{0}}$ on $U$. Then  $H_{0} := \{x \in U \ / \ f_{0}(x) = 0 \}$}. Then there exists near $(0,p)$ in $D \times M$ a smooth hypersurface $H$ such that for each $s$ near enough $0$ the fiber $X_{s}$ is also contained in $H$ in a fixed open set around $(0,p)$.
\end{lemma}

\parag{Proof} Let  $f_{0}$ be  a holomorphic function in an open neighbourghood $U$ of $p$ such that $(df_{0})_{p} \not= 0$ and $f_{0}$ is a section of $\I_{X_{0}}$ on $U$ and define $H_{0} := \{x \in U \ / \ f_{0}(x) = 0 \}$. Then for $0 \in D' \subset D$ small enough there exist  $f \in \Gamma(D'\times U, \I_{X})$ inducing $f_{0}$ on $U$. Then we have $df_{0,p}\not= 0$ and $H := \{(s,x) \in (D \times M)\cap(D'\times U)\}$ is a smooth hypersurface near $(0,p)$ which contains $X$ so $X_{s} \cap (D'\times U)$ for $s \in D'$.$\hfill \square$

\parag{Application} Assume that in the situation of the previous lemma we know that $red(X_{0})$, the reduced fiber of $\pi$ at $0$,  is contained near $p$ in a smooth hypersurface in $M$\footnote{this is equivalent to the fact that its Zariski tangent space at $p$ has codimension $1$ in $T_{M,p}$.}, and that $X$ is not contained in a smooth hypersurface in $D\times M$ near $(0,p)$. Then the fiber $X_{0}$ cannot be reduced.

\parag{Example} Let $M := \C^{3}$ and define
 $$X := \{\big(s,(x,y,z) \big) \in D \times \C^{3} \ / \   x = y = 0 \quad {\rm or} \quad  x - s = z = 0 \} $$
  with its reduced structure in $D \times M$. Then $X$ is not contained in a smooth hypersurface near $(0,(0,0,0))$ because if $\varphi$ is an holomorphic function near $(0,(0,0,0))$ vanishing on $X$ we have $\varphi = a.x + b.y = c.(x-s) + d.z$ with holomorphic functions $a,b,c,d$ and as $x,y,z,s$ is a regular sequence, it implies that $\varphi$ is in the square of the maximal ideal at the origin. This implies that the Zariski tangent space to $X$ at the origin is $\C^{4}$. But remark that $red(X_{0})$ is contained in $\{ x = 0 \}$. So we conclude that the fiber $X_{0}$ is not reduced. A direct (and easy) calculation gives that $x$ is a non trivial  nilpotent element in sheaf  $\mathcal{O}_{X_{0}} \simeq \mathcal{O}_{\C^{3}}\big/(x^{2}, x.y, x.z, y.z)$.\\
  
  Note that the example above describe exactly the situation that happens if we try to deform the example of non K\"ahler manifold constructed in the begining of this section by moving the two branches of the curve $C$ near the origin in order that for $s \not= 0$ they dont meet. So the fiber at $s = 0$ of such a family will not be reduced ! 

\parag{Remark} For any pair of smooth families of smooth curves in $\C^{3}$ parametrized by $D$ such that the graph $\Gamma$ and $\Gamma'$ are smooth surfaces in $D \times \C^{3}$ which meets  transversally in a point $(0,p)$ the same phenomenon occurs because in $\C^{4}$ the union of two transversal $2-$planes at $0$ (with its reduced structure) is not locally contained in a smooth hypersurface near $0$\footnote{Two surfaces in a smooth $3-$fold which meet, meet in codimension at most equal to $2$.\\
But remark that in the case we are interested with, they are locally contained in a (singular)  hypersurface like $x.u = y.v$ where $(x,y,u,v)$ are local coordinates near the origin.} and the Zariski tangent space of two transversal curves in $\C^{3}$ has always codimension $1$ in $\C^{3}$. So the union of two transversal curves at $0$ in $\C^{3}$ is always locally contained in a smooth (hyper-)surface.

\parag{Three curves example} Consider the analytic set
 $$Y := \{x = y = 0\} \cup \{y = z = 0\} \cup \{x - s = z = 0\} $$
  in $\C^{4}$ with its reduced structure. We have $\I_{Y} = ( y.(x-s), x.z, y.z )$ and the ideal $\I_{Y_{0}} = (x.y, x.z, y.z)$ is clearly the reduced ideal of the fiber $Y_{0}$ for the projection on $D$ corresponding to the coordinate $s$.\\
We want to prove that the map $\pi : Y \to D$ given by the $s-$projection is flat and has a {\bf reduced} fiber at $s = 0$. So let us prove that $\mathcal{O}_{Y}$ has no $s-$torsion. Consider a holomorphic function $\varphi$ near $0$ in $\C^{4}$ such that $s.\varphi$ is in $\I_{Y}$ near $0$. Then modulo $(z)$ we have $s.[ \varphi] \in [y.(x - s)]$, and, as $s, y(x - s)$ is a regular sequence in $\C^{3}$ we obtain that $\varphi = a.y.(x - s) + z.b $ where $a, b$ are holomorphic. Now we have $s.z.b \in \I_{Y}$ and this gives modulo $(y)$ that $s.[z.b] = c.[x.z] $. But $z$ does not divide $0$ modulo $(y)$ so $s.[b]= c.[x]$ and then, as $s, x,y$ is a regular sequence  $b = u.x + v.y$
and $\varphi = a.y.(x - s) + u.z.x + v.z.y \in \I_{Y} $. So $Y$ is flat on $D$ and its fiber at $s =0$ is reduced because its ideal is generated by $x.y, x.z, y.z$.\\

This enlight why H. Hironaka considers $3$ transversal curves at the point $Q$ and let one of the $3$ curves moving outside this point; this allows to preserve in such a family the construction given at $s = 0$. Of course, this does not replace the proof of the section $4$ in  [H.61] but explains that with three curves it is possible, although it does not work with two.

\section{A last result}
Choose now in $\mathbb{P}_{3}$ two smooth connected  curves $C_{1}$ and $C_{2}$ meeting transversally at two points $P$ and $Q$ and perform for each point the previous construction consisting in the blowing-up of one branch after the other in a small ball around each point, but reversing the order : at $P$ we blow-up first $C_{1}$ and after the strict transform of $C_{2}$ but near $Q$ we blow-up first $C_{2}$ and then the strict transform of $C_{1}$. Of course, outside the two balls we just blow-up $C_{1}$ and $C_{2}$ which are disjoint. So the all thing patch in a complex manifold $Y$ and in a modification $\pi : Y \to \mathbb{P}_{3}$. We shall denote by $L_{1}$ the generic fiber of $E_{1} :=\pi^{-1}(C_{1})$ on $C_{1}$ and $L_{2}$ the generic fiber of $E_{2} := \pi^{-1}(C_{2})$ on $C_{2}$.\\
Then we find over $P$ two smooth rational  curves $A$ and $B$ and we have $L_{1} \sim A + B $ and $L_{2}\sim B$. In a similar way, we find over $Q$ two smooth rational curves $C$ and $D$ which satisfy $ L_{2} \sim C + D $ and $L_{1} \sim D$. This gives $A + C \sim 0$ in $Y$. So we have again a non k{\"a}hler smooth Moishezon manifold. See picture 6.\\

The main interest in this example is that we have two disjoint curves, such that the blow-up of one of the two gives back a projective manifold.

\begin{lemma}\label{dernier} 
If we blow-up $A$ {\bf or}  $C$ in the manifold $Y$ defined above we obtain a projective manifold.
\end{lemma}

\parag{Proof} In fact, as $A$ is the pull-back of the point $P$ in the first blow-up, to blow-up $Y$ along $A$ will give the same result than to begin by the blow-up of $P$ and then to blow-up the strict transform of $C_{1}$ and then the strict transform of $C_{2}$. But after blowing-up $P$ the strict transforms of $C_{1}$ and $C_{2}$ are disjoint in the pull-back of a small ball around $P$  and the order does not matter then. So after the blow-up of $P$ we can simply blow-up in a global manner the strict transform of $C_{2}$ and then the strict transform of $C_{1}$ respecting the order choosen in the ball around $Q$. So we have performed globally $3$ blow-up of smooth projective subvarieties and the result is projective. The situation is analog for $C$.$\hfill \square$\\

This example shows that, in some sense, the non k{\"a}hlerianity of this manifold is ``concentrated around $A$'' or ``around $C$''  but $A$ and $C$ are disjoint ... 

\bigskip

\bigskip

{\bf Bibliography}: \begin{itemize}
\item{[H.61]} Hironaka, H. {\it An example of a non-k{\"a}hlerian complex-analytic deformation of k{\"a}hlerian complex structures} Annals of Math. 75, 1, (1962) p.190-208.
\end{itemize}

\begin{figure}[h]
\section{Pictures} Remember that these pictures try to represent subsets in a $3-$dimensional complex space. So they are necessarily ``false'' and can only help to understand using imagination.

\bigskip

\bigskip

\begin{center}
\includegraphics[width=18cm]{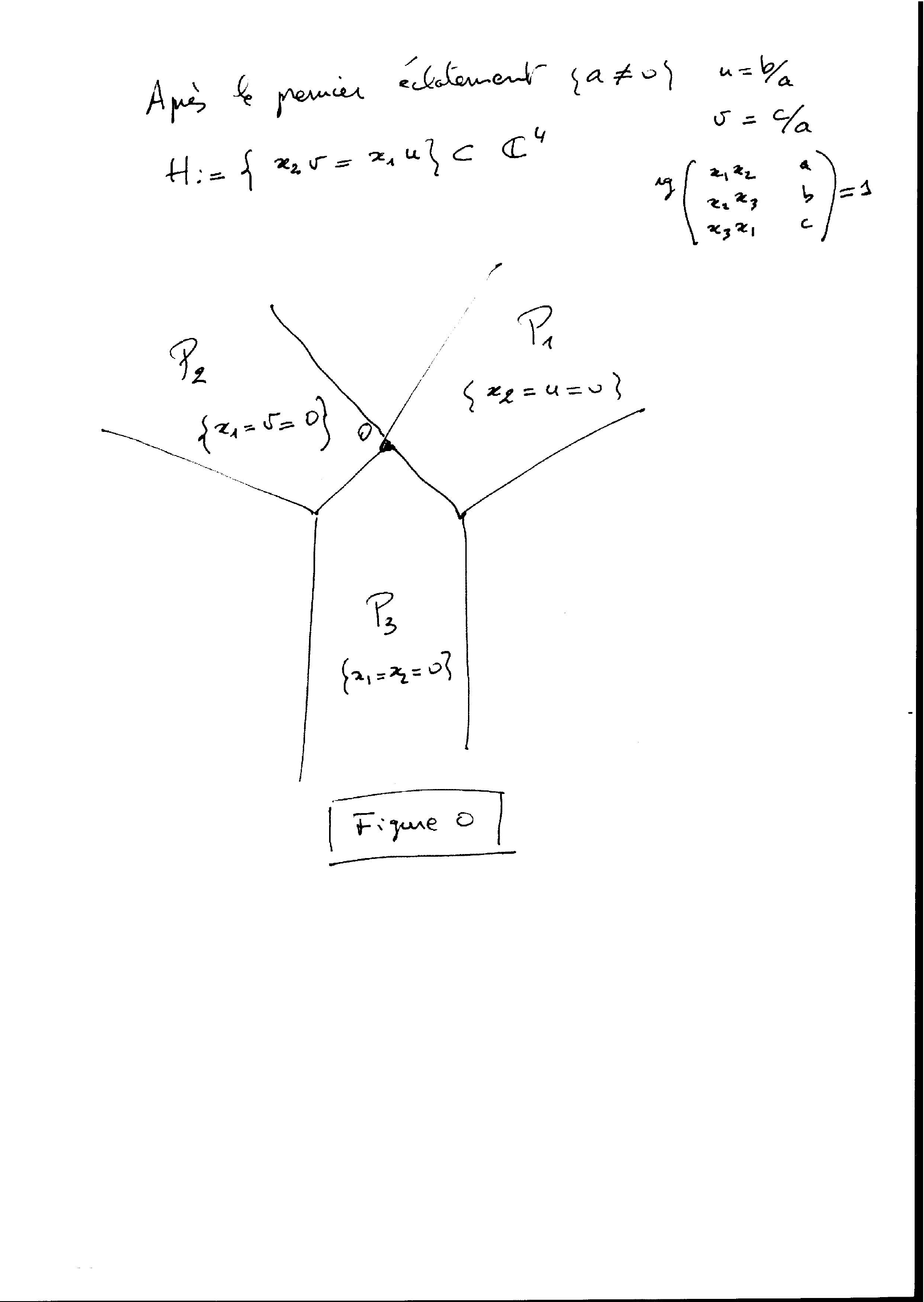}\\
\end{center}
\end{figure}

\bigskip

\begin{figure}[h]
\begin{center}
\includegraphics[width=20cm]{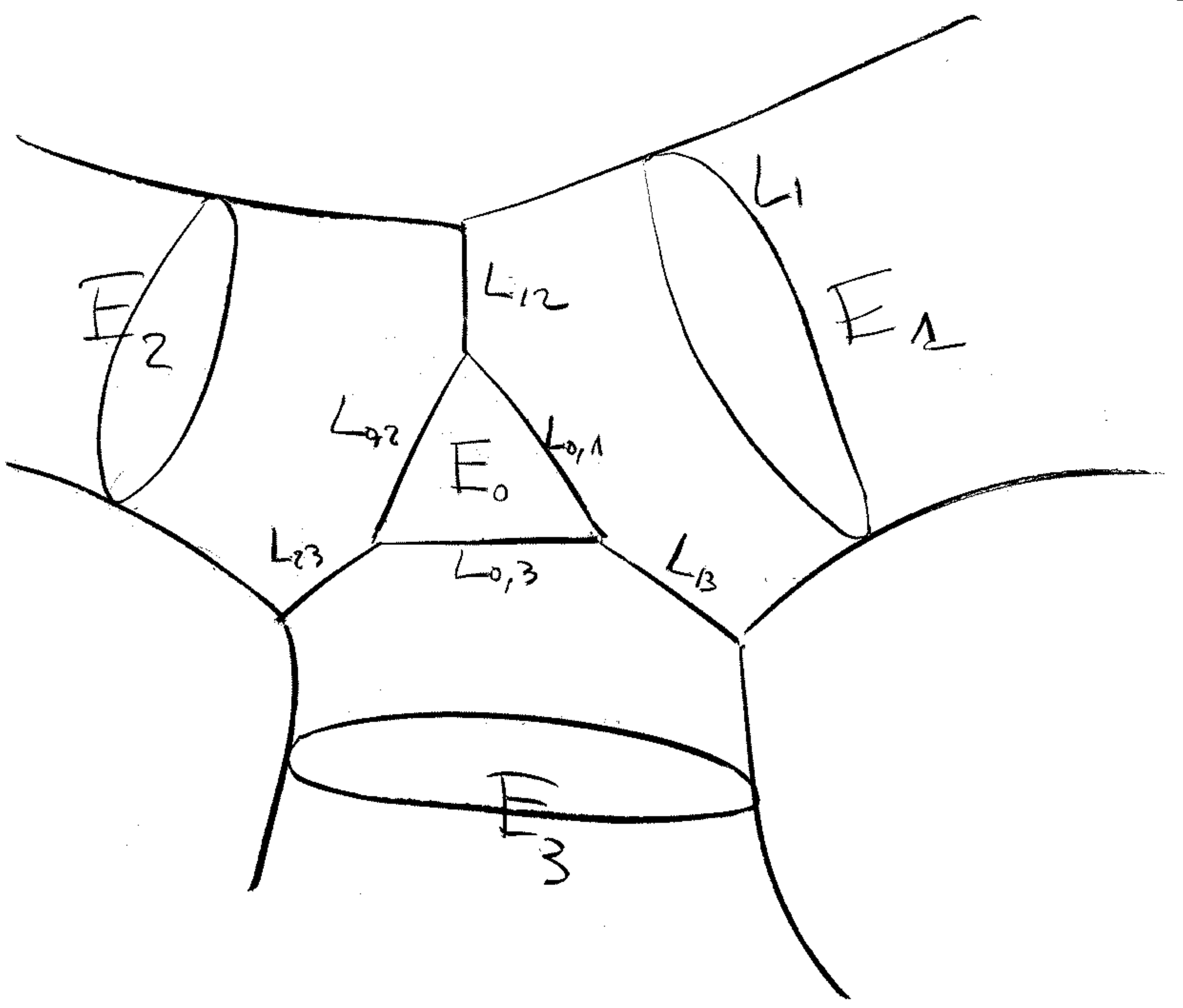}\\
Figure 1.
\begin{align*}
& L_{1} \sim L_{1,2} + L_{0,1} + L_{1,3}  \\
& L_{2} \sim L_{2,3} + L_{1,2} + L_{0,2}  \\
& L_{3} \sim L_{3,1} + L_{2,3} + L_{0,3}
\end{align*}
\end{center}
\end{figure}

\newpage

\begin{figure}
\includegraphics[width=20cm]{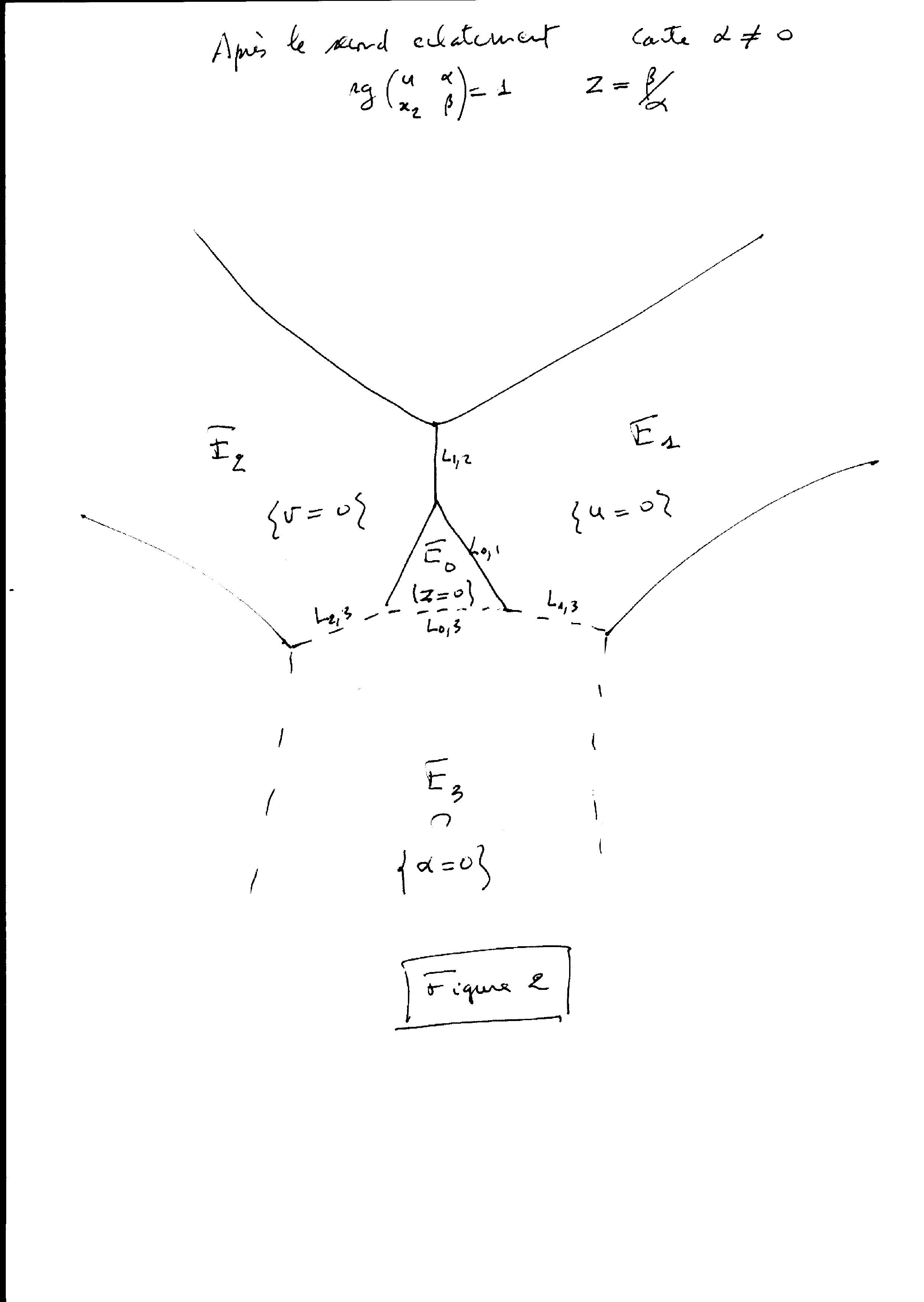}\\
\end{figure}

\begin{figure}
\begin{center}
\includegraphics[width=18cm]{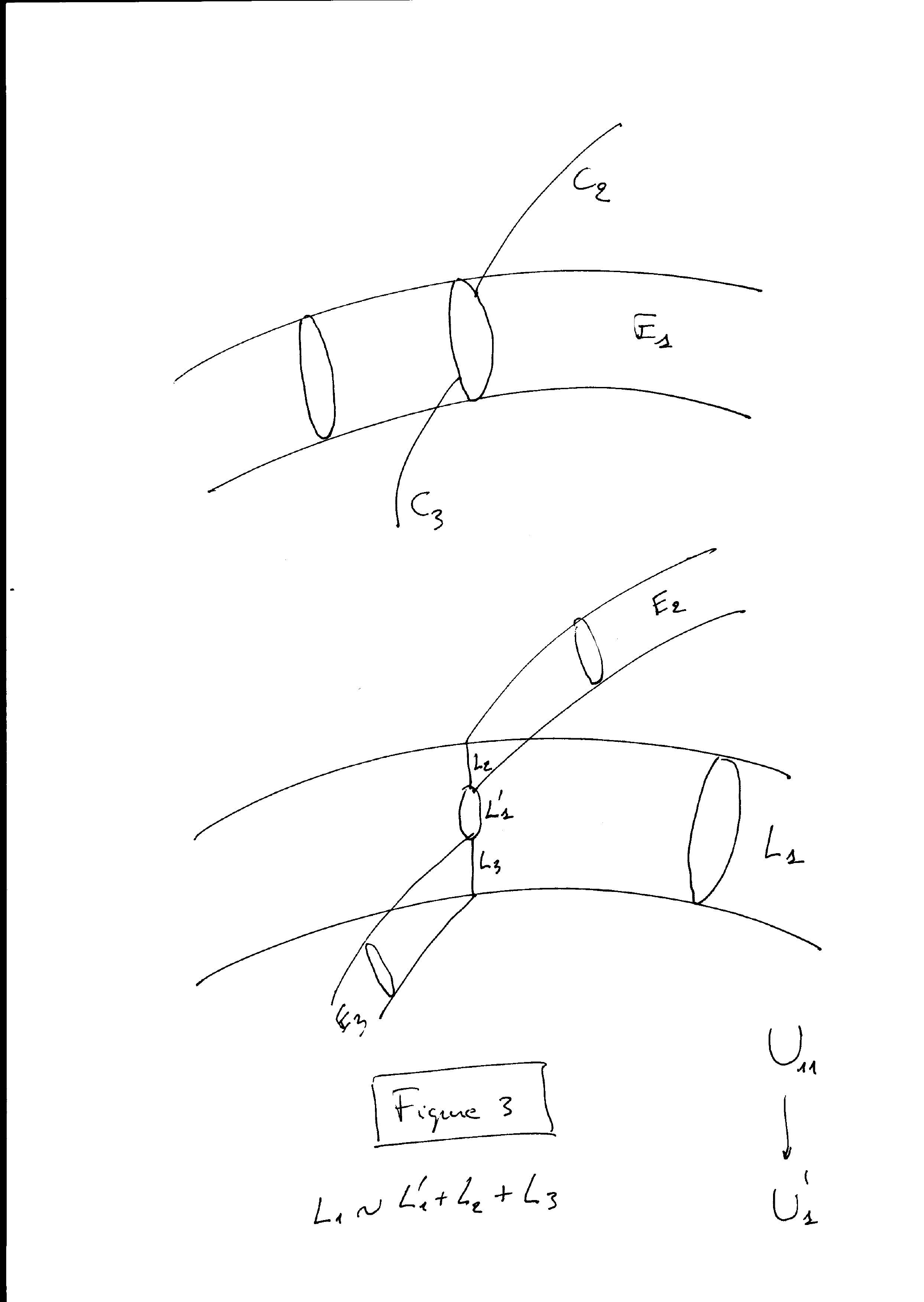} \\
\end{center}
\end{figure}

\begin{figure}
\begin{center}
\includegraphics[width=22cm]{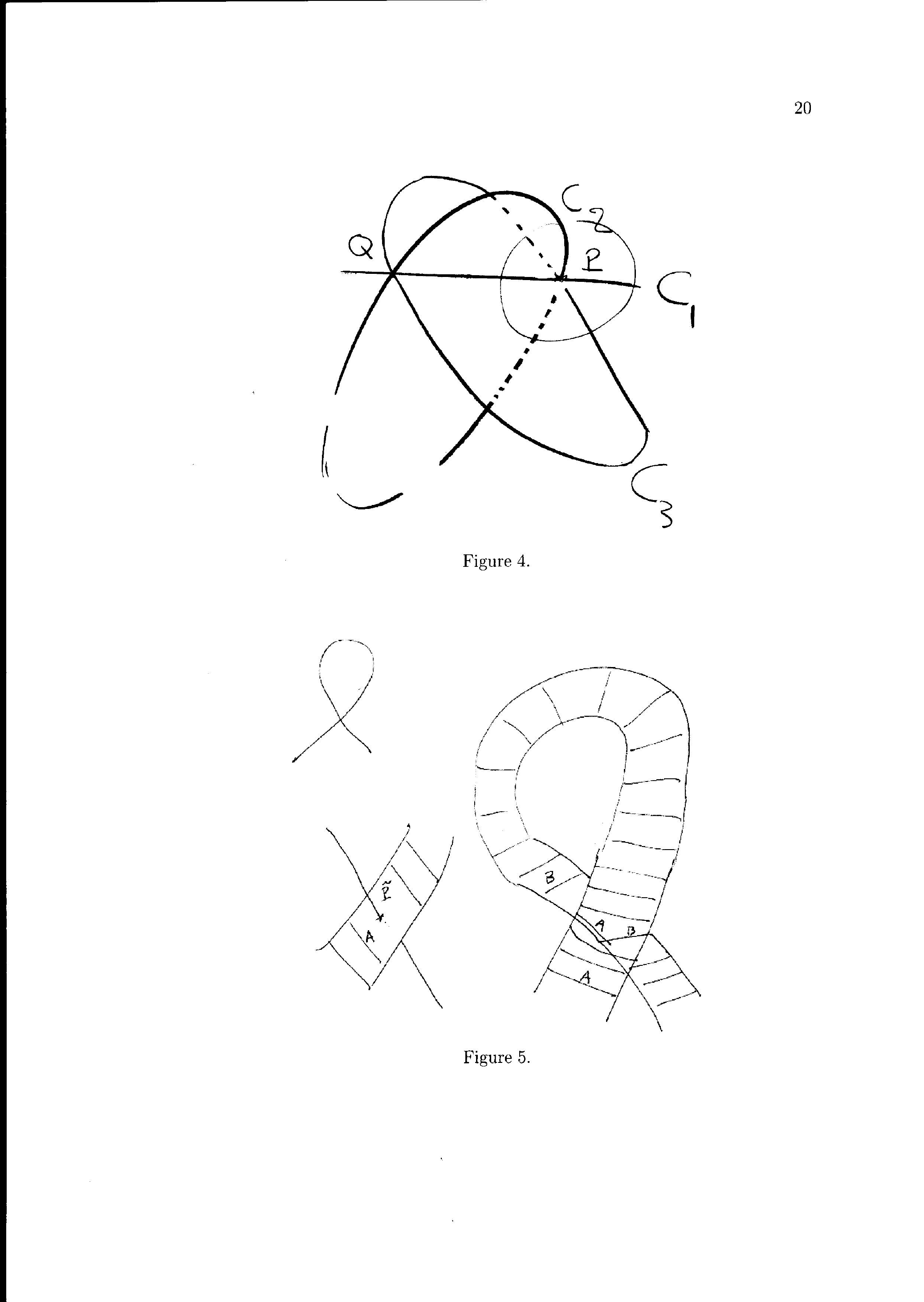}\\
\end{center}
\end{figure}


\begin{figure}
\begin{center}
\includegraphics[width=18cm]{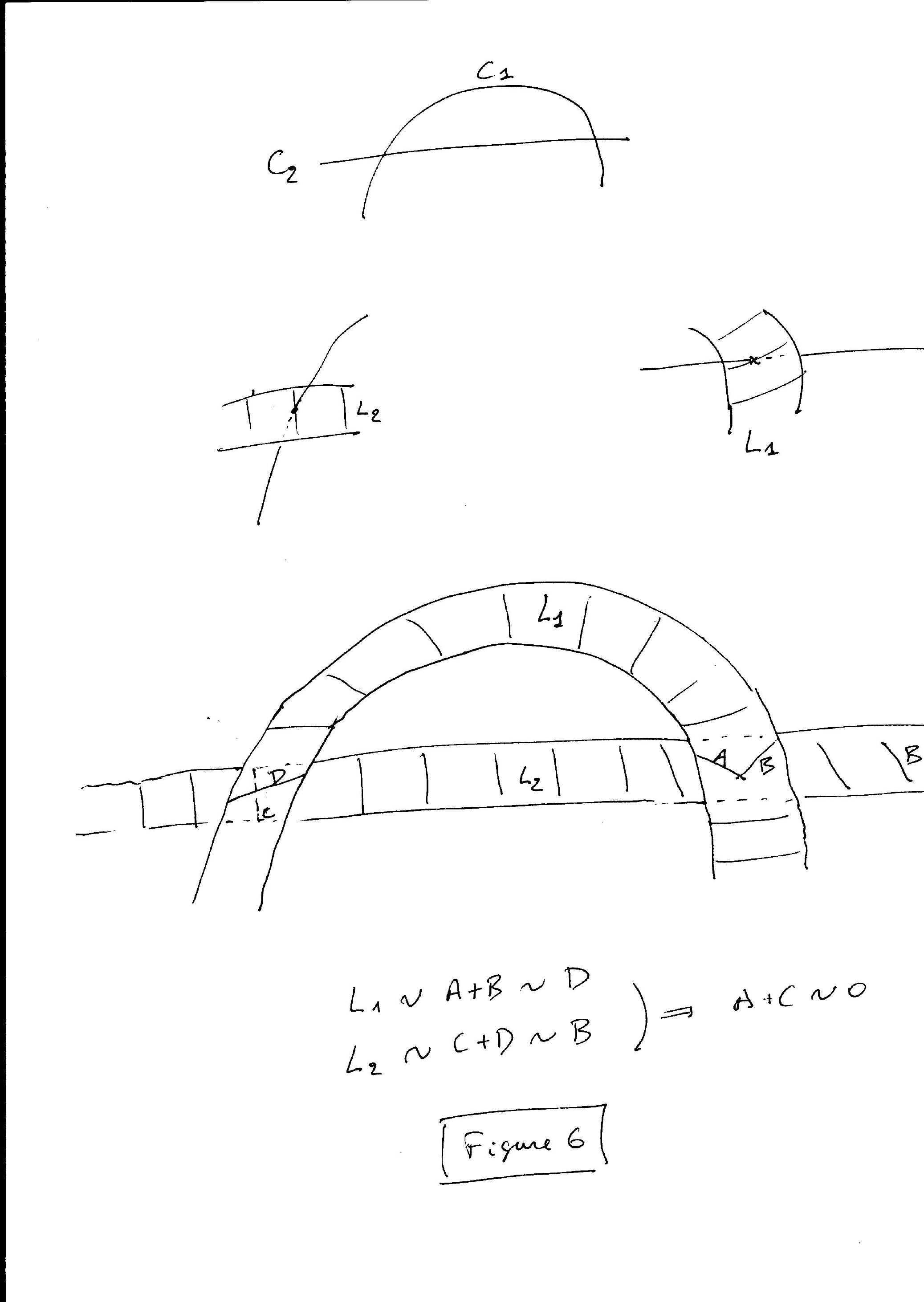}\\
\end{center}
\end{figure}

\end{document}